\documentclass[a4paper,12pt]{amsart}
\usepackage{amssymb}
\usepackage{ifthen}
\usepackage[dvips]{graphicx}
\nonstopmode \numberwithin{equation}{section}
\usepackage{amssymb}
\usepackage{ifthen}
\usepackage{graphicx}
\usepackage{amsmath}
\usepackage[T1]{fontenc} 
\usepackage[utf8]{inputenc}
\usepackage[usenames,dvipsnames]{color}
\usepackage{color}
\usepackage[english]{babel}
\usepackage{fancyhdr}
\usepackage{fancybox}
\usepackage{tikz}

\nonstopmode \numberwithin{equation}{section}
\theoremstyle{plain}
\newtheorem{thm}{Theorem}[section]
\newtheorem{cor}[equation]{Corollary}
\newtheorem{lem}[equation]{Lemma}
\newtheorem{prop}{Proposition}

\newtheorem{conj}{Conjecture}

\theoremstyle{definition}
\newtheorem{defn}{Definition}[section]

\newtheorem{prob}{Problem}
\newtheorem{rem}{Remark}[section]

\newcounter{minutes}
\divide\time by 60
\newcounter{hours}
\multiply\time by 60
\addtocounter{minutes}{-\time}

\newcounter {own}
\def\theown {\thesection       .\arabic{own}}

\newenvironment{pf}[1][]{%
	\vskip 3mm
	\noindent
	\ifthenelse{\equal{#1}{}}%
	{{\slshape Proof. }}%
	{{\slshape #1.} }%
}%
{\qed\bigskip}

\newcounter{alphabet}

\def\be{\begin{equation}}
	\def\ee{\end{equation}}

\newcommand{\bee}{\begin{enumerate}}
	\newcommand{\eee}{\end{enumerate}}

\newcommand{\blem}{\begin{lem}}
	\newcommand{\elem}{\end{lem}}
\newcommand{\bthm}{\begin{thm}}
	\newcommand{\ethm}{\end{thm}}
\newcommand{\bcor}{\begin{cor}}
	\newcommand{\ecor}{\end{cor}}
\newcommand{\beg}{\begin{examp}}
	\newcommand{\eeg}{\end{examp}}
\newcommand{\begs}{\begin{examples}}
	\newcommand{\eegs}{\end{examples}}

\newcommand{\bdefn}{\begin{defn}}
	\newcommand{\edefn}{\end{defn}}

\newcommand{\bprob}{\begin{prob}}
	\newcommand{\eprob}{\end{prob}}
\newcommand{\bei}{\begin{itemize}}
	\newcommand{\eei}{\end{itemize}}

\newcommand{\bcon}{\begin{conj}}
	\newcommand{\econ}{\end{conj}}
\newcommand{\bcons}{\begin{conjs}}
	\newcommand{\econs}{\end{conjs}}
\newcommand{\bprop}{\begin{prop}}
	\newcommand{\eprop}{\end{prop}}
\newcommand{\br}{\begin{rem}}
	\newcommand{\er}{\end{rem}}
\newcommand{\brs}{\begin{rems}}
	\newcommand{\ers}{\end{rems}}
\newcommand{\bo}{\begin{obser}}
	\newcommand{\eo}{\end{obser}}
\newcommand{\bos}{\begin{obsers}}
	\newcommand{\eos}{\end{obsers}}
\newcommand{\bpf}{\begin{pf}}
	\newcommand{\epf}{\end{pf}}
\newcommand{\ba}{\begin{array}}
	\newcommand{\ea}{\end{array}}
\newcommand{\beq}{\begin{eqnarray}}
	\newcommand{\beqq}{\begin{eqnarray*}}
		\newcommand{\eeq}{\end{eqnarray}}
	\newcommand{\eeqq}{\end{eqnarray*}}

\begin{document}
	
	\title{Improved Bohr radius for the class of starlike log-harmonic mappings}
	
	\author{Molla Basir Ahamed}
	\address{Molla Basir Ahamed,
		School of Basic Science,
		Indian Institute of Technology Bhubaneswar,
		Bhubaneswar-752050, Odisha, India.}
	\email{mba15@iitbbs.ac.in}
	
	\author{Vasudevarao Allu}
	\address{Vasudevarao Allu,
		School of Basic Science,
		Indian Institute of Technology Bhubaneswar,
		Bhubaneswar-752050, Odisha, India.}
	\email{avrao@iitbbs.ac.in}

	\subjclass[{AMS} Subject Classification:]{Primary 30C45, 30C50, 30C80}
	
	\def\thefootnote{}
	\footnotetext{ {\tiny File:~\jobname.tex,
			printed: \number\year-\number\month-\number\day,
			\thehours.\ifnum\theminutes<10{0}\fi\theminutes }
	} \makeatletter\def\thefootnote{\@arabic\c@footnote}\makeatother

\begin{abstract}
Let $ \mathcal{H}(\mathbb{D}) $ be the linear space of analytic functions on the unit disk $ \mathbb{D}=\{z\in\mathbb{C}: |z|<1\} $ and let $ \mathcal{B}=\{w\in \mathcal{H}(\mathbb{D}: |w(z)|<1)\} $. The classical Bohr's inequality states that if a power series $ f(z)=\sum_{n=0}^{\infty}a_nz^n $ converges in $ \mathbb{D} $ and $ |f(z)|<1 $ for $ z\in\mathbb{D} $, then 
\begin{equation*}
	\sum_{n=0}^{\infty}|a_n|r^n\leq 1\;\;\mbox{for}\;\; r\leq \frac{1}{3}
\end{equation*}
and the constant $ 1/3 $ is the best possible.  The constant $ 1/3 $ is known as Bohr radius. A function $ f : \mathbb{D}\rightarrow\mathbb{C} $ is said to be log-harmonic if there is a $ w\in\mathcal{B} $ such that $ f $ is a non-constant solution of the non-linear elliptic partial differential equation
\begin{equation*}
	\bar{f}_{\bar{z}}(z)/\bar{f}(z)=w(z)f_{z}(z)/f(z).
\end{equation*}
The class of log-harmonic mappings is denoted by $ \mathcal{S}_{LH} $. The set of all starlike log-harmonic mapping is defined by
\begin{equation*}
\mathcal{ST}_{LH}=\bigg\{f\in\mathcal{S}_{LH}:\frac{\partial}{\partial\theta}{\rm Arg}(f(e^{i\theta}))={\rm Re}\left(\frac{zf_{z}-\bar{z}f_{\bar{z}}}{f}\right)>0\;\; \mbox{in}\;\; \mathbb{D}\bigg\}.
\end{equation*}
In this paper, we study several improved Bohr radius for the class $ \mathcal{ST}^{0}_{LH} $, a subclass of $ \mathcal{ST}_{LH} $, consisting of functions $ f\in\mathcal{ST}_{LH} $ which map the unit disk $ \mathbb{D} $ onto a starlike domain (with respect to the origin).
\end{abstract}

\maketitle
\pagestyle{myheadings}
\markboth{Molla Basir Ahamed and Vasudevarao Allu}{Improved Bohr radius for the class of starlike log-harmonic mappings}

\section{Introduction and preliminaries}
 \par A complex-valued function $ f $ in $ \mathbb{D} $ is said to be harmonic if it satisfies the Laplace equation $ \Delta f=4f_{z\bar{z}}=0 $ in $ \mathbb{D} $. Every harmonic function $ f $ in $ \mathbb{D} $ has the unique canocical form $ f=h+\bar{g} $, where $ h $ and $ g $ are analytic in $ \mathbb{D} $ with $ g(0)=0 $. Every analyitc function is a harmonic function. Let $ \mathcal{H} $ be the class of all complex-valued harmonic functions $ f=h+\bar{g} $ defined on $ \mathbb{D} $, where $ h $ and $ g $ are analytic in $ \mathbb{D} $ with the normalization $ h(0)=h^{\prime}(0)-1=0 $ and $ g(0)=0. $ Here $ h $ is called analytic part and $ g $ is called co-analytic part of $ f $.\vspace{2mm}

Harmonic mappings play the natural role in parameterizing minimal surfaces in the context of differential geometry. Planner harmonic mappings have application not only in the differential geometry but also in various field of engineering, physics, operations research and other intriguing aspects of applied mathematics. The theory of harmonic functions has been used to study and solve fluid flow problems \cite{aleman-2012}. The theory of univalent harmonic functions having prominent geometric properties like starlikeness, convexity and close-to-convexity appear naturally while dealing with planner fluid dynamical problems. For instance, the fluid flow problem on a convex domain satisfying an interesting geometric property has been extensively studied by Aleman and Constantin \cite{aleman-2012}. With the help of geometric properties of harmonic mappings, Constantin and Martin \cite{constantin-2017} have obtained a complete solution of classifying all two dimensional fluid flows.\vspace{2mm}

Let $ \mathcal{H}(\mathbb{D}) $ be the class of analytic functions $ f(z)=\sum_{n=0}^{\infty}a_nz^n $ in the unit disk $ \mathbb{D}$. The origin of the Bohr phenomenon lies in the seminal work by Harald Bohr \cite{Bohr-PLMS-1914}, which include the following result: Let $ f\in\mathcal{H}(\mathbb{D}) $ satisfies $ |f(z)|<1 $ for all $ z\in\mathbb{D} $, then $ \sum_{n=0}^{\infty}|a_n||z|^n\leq 1 $ for all $ z\in\mathbb{D} $ with $ |z|=r\leq 1/3 $ and the constant $ 1/3 $ is the best possible. For $ f\in\mathcal{H}(\mathbb{D}) $, the majorant series is denoted by $ M_f(r) $ and is defined by $ M_f(r)=\sum_{n=0}^{\infty}|a_n||z|^n $.  Bohr actually obtained the inequality $ M_f(r)\leq 1 $ for $ |z|\leq 1/6 $, but subsequently later, M. Riesz, I. Schur and F. Weiner, independently established this inequality for $ |z|\leq 1/3 $ and the constant $ 1/3 $ cannot be improved \cite{Paulson-Popescu-Singh-2002}. The constant $ r_0=1/3 $ is called the Bohr radius and the inequality $ M_f(r)\leq 1 $ is called Bohr inequality for bounded analytic functions in the unit disk $ \mathbb{D}. $ Moreover, for the function $ \phi_a $ defined by
\begin{equation*}
	\phi_a(z)=\frac{a-z}{1-az},\;\; a\in [0,1)
\end{equation*}
it follows that $ M_{\phi_a}(r)>1 $ if, and only if, $ r>1/(1+2a) $, for which $ a\rightarrow 1 $ shows that $ 1/3 $ is optimal.\vspace{2mm}

Using the Euclidian distance $ d $, the Bohr inequality for $ f\in\mathcal{H}(\mathbb{D}) $ can be written as 
\begin{equation}\label{e-1.1}
	d\left(\sum_{n=0}^{\infty}|a_nz^n|, |a_0|\right)=\sum_{n=1}^{\infty}|a_nz^n|\leq 1-|f(0)|=d(f(0),\partial\mathbb{D}),
\end{equation}
where $ \partial\mathbb{D} $ is the boundary of the unit disk $ \mathbb{D} $.\vspace{2mm}

\noindent Let $ \mathcal{M} $ be a class of analytic functions $ f(z)=\sum_{n=0}^{\infty}a_nz^n $ which map the unit disk $ \mathbb{D} $ into a domain $ \Omega\subset\mathbb{C} $. We say the class $ \mathcal{M} $ satisfies Bohr phenomenon if there exists $ r^{*} $ such that \eqref{e-1.1} holds for $ |z|=r\leq r^{*} $. The largest such $ r^{*} $ is called the Bohr radius for the class $ \mathcal{M} $. \vspace{2mm}

In the recent years, studying Bohr inequalities  become an interesting topic of research for the functions of one as well as several complex variables. The notion of Bohr inequality has been generalized to several complex variables (see \cite{Aizenberg-PAMS-2000,Aizenberg-Aytuna-JMAA-2001,Boas-Khavinshon-PAMS-1997,Liu-Ponnusamy-PAMs-2020}), to planner harmonic mappings (see \cite{Evdoiridis-Ponnusamy-IM-2019,Kayumov-Ponnusamy-MN-2018,Kayumov-Ponnusamy-JMAA-2018}) to polynomials (see \cite{Fournier-JMAA-2008}), to the solutions of elliptic partial differential equations (see \cite{Abdulhadi-Hengartner-1989,Ali-Abdulhadi-Ng-CVEE-2016}), to elliptic equations (see \cite{Aizenberg-Tarkhanov-PLMS-2001}), to vector valued functions and operaotor valued functions (see \cite{Bhowmik-Das-2019,Bhowmik-Das-arxive-2020}), to analytic functions in norm linear spaces (see \cite{Dahlnere-Khavinshon-CMFT-2004}) and in a more abstract setting (see \cite{Aizenberg-Aytuna-Djakov-PAMS-2000}). In $ 1977 $, Boas and Khavinshon \cite{Boas-Khavinshon-PAMS-1997} extended the Bohr inequality to several complex varibales by finding multidimensional Bohr radius. Bohr's theorem attracted a greater interest after it was used by Dixon \cite{Dixon-BLMS-1995} in $ 1995 $ to characterize Banach algebras that satisfy von Neumann inequality. The generalization of Bohr's theorem become now-a-days an active topic of research. In $ 2001 $, Aizenberg \emph{et al.} \cite{Aizenberg-Aytuna-JMAA-2001}, and in $ 2013 $ Aytuna and Djakov \cite{Aytuna-Djkov-BLMS-2013} studied the Bohr property of holomorphic functions while Paulsen \emph{et al.}
 \cite{Paulson-Popescu-Singh-2002} extended the Bohr inequality to Banach algebra. The relevance between Banach theory and Bohr's theorem was explored in \cite{Blasco-2009,Defant-2003,Dixon-BLMS-1995}.
 \vspace{2mm}
 
 Recently, Ali and Ng \cite{Ali-Ng-CVEE-2018} have extended the classical Bohr inequality in the Poincare disk model of hyperbolic plane. Kayumove and Ponnusamy \cite{Kayumov-Ponnusamy-JMAA-2018} have determined the Bohr radius for the class of analytic functions $ f(z)=z^m\sum_{k=0}^{\infty}a_{kp}z^{kp} $, $ p\geq m\geq 0 $ with $ |f(z)|\leq 1 $. In $ 2018 $, Kayumov \emph{et al.} \cite{Kayumov-Ponnusamy-MN-2018} introduced the idea of $ p $- Bohr radius for harmonic functions and obtained the $ p $- Bohr radius for the class of odd harmonic functions. Kayumov \emph{et al.} \cite{Kayumov-Ponnusamy-MN-2018} have obtained the Bohr radius for the class of analytic Bloch functions and harmonic functions. Alkhaleefah \emph{et al.} \cite{Alkhaleefah-Kayumov-Ponnusamy-PAMS-2019} have studied the Bohr radius for the class of quasi-subordinate functions which in particular gives the classical Bohr radius. Number of improved versions of the classical Bohr inequality have been proved in \cite{Kayumov-Ponnusamy-CAMS-2020}.\vspace{2mm}
 
 \par We now define Bohr radius in subordination and bounded harmonic classes.  Let $ f $ and $ g $ be two analytic functions in the unit disk $ \mathbb{D} $. We say that $ g $ is subordinate to $ f $ if there exists an anlytic function $ \phi : \mathbb{D}\rightarrow\mathbb{D} $ with $ \phi(0)=0 $ so that $ g=f\circ \phi $ and it is denoted by $ f\prec g $. If $ g $ is univalent and $ f(0)=g(0) $ then $ f(\mathbb{D})\subset g(\mathbb{D}) $. We denote the class of all functions subordinate to a fixed function $ f $ by $ \mathcal{S}(f) $ and $ f(\mathbb{D})=\Omega $. The class $ \mathcal{S}(f) $ is said to have Bohr's phenomenon if for any $ g(z)=\sum_{n=0}^{\infty}b_nz^n\in \mathcal{S}(f) $ and $ f(z)=\sum_{n=0}^{\infty}a_nz^n $ there is a $ r_0 $ in $ (0,1] $ such that
 \begin{equation}\label{e-2.1}
 	\sum_{n=0}^{\infty}|b_nz^n|\leq d(f(0),\partial\Omega)\;\;\mbox{for}\;\; |z|<r_0. 
 \end{equation}
 \noindent In $ 2010 $, it was established by Abu-Muhanna \cite[Theorem]{Muhanna-CVEE-2010} that the class $ \mathcal{S}(f) $ has Bohr phenomenon when $ f $ is univalent in $ \mathbb{D} $. In particular, the following interesting result was obtained.
 \begin{thm}\cite{Muhanna-CVEE-2010}\label{th-2.1}
 	If $ g(z)=\sum_{n=0}^{\infty}b_nz^n\in \mathcal{S}(f) $ and $ f(z)=\sum_{n=0}^{\infty}a_nz^n $ is univalent, then\begin{equation}\label{e-2.3}
 		\sum_{n=1}^{\infty}|b_nz^n|\leq d(f(0),\partial\Omega)\;\; \mbox{for}\;\; |z|\leq r_0=3-\sqrt{8}=0.17157.
 	\end{equation}
 	Here $ r_0 $ is sharp for the Koebe function $ f_K(z)=z/(1-z)^2. $
 \end{thm}
 \noindent In \cite{Muhanna-CVEE-2010}, Abu-Muhanna has proved the following lemma to find the lower bound of the distance $ d(f(0),\partial\Omega) $.
 \begin{lem}\cite{Muhanna-CVEE-2010}\label{lem-2.3}
 	Let $ f(z)=\sum_{n=0}^{\infty}a_nz^n $ be an analytic univalent function from $ \mathbb{D} $ onto a simply connected domain $ \Omega $. Then 
 	\begin{equation}\label{e-2.4}
 		\frac{1}{4}|f^{\prime}(0)|\leq d(f(0),\partial\Omega)\leq |f^{\prime}(0)|.
 	\end{equation}
 \end{lem}

Next we discuss improved Bohr radius for starlike log-harmonic mappings. A function $ f : \mathbb{D}\rightarrow\mathbb{C} $ is said to be log-harmonic if there is a $ w\in\mathcal{B} $ such that $ f $ is a non-constant solution of the non-linear elliptic partial differential equation
\begin{equation}\label{e-2.6}
	\bar{f}_{\bar{z}}(z)/\bar{f}(z)=w(z)f_{z}(z)/f(z),
\end{equation} 
where the second dilation function $ w $ is such that $ |w(z)|<1 $ for all $ z\in\mathbb{D} $. The Jacobian
\begin{equation*}
	J_f=|f_{z}|^2-|f_{\bar{z}}|^2=|f_z|^2(1-|w(z)|^2)
\end{equation*}
is positive, and therefore all the non-constant log-harmonic mappings are sense-preserving and open in $ \mathbb{D} $.\vspace{2mm}

\par In $ 2013 $, Li \textit{et al.} \cite{Li-Ponnusamy-Wang-BMMS-2013} proved a necessary and sufficient condition for a function to be log-$ p $-harmonic and also studied local log-$ p $-harmonic mappings. Mao \textit{et al.} \cite{Mao-Ponnusamy-Wang-CVEE-2013} have established Schwarz' lemma for log-harmonic mappings, through which they proved two versions of the Landau's theorem for these functions. In $ 2018 $, Liu and Ponnusamy \cite{Liu-Ponnusamy-2020} obtained the coefficient estimates and hence studied Bohr radius for log-harmonic mappings. Inner mapping radius by constructing a family of $ 1 $-slit log-harmonic mappings have been established in \cite{Liu-Ponnusamy-2020}. Several interesting properties have been established in \cite{Liu-Ponnusamy-2020} of log-harmonic mappings. In $ 2019 $, Liu and Ponnusamy \cite{Liu-Ponnuamy-arxive-2019} obtained the precise ranges of log-harmonic Koebe mapping, log-harmonic right half-plane mapping and log-harmonic two-slits mappings. Further, the coefficient estimates for univalent log-harmonic starlike mappings has been established in \cite{Liu-Ponnuamy-arxive-2019}.\vspace{2mm}

\noindent Let $ h_0 $ and $ g_0 $ be two functions defined by 
\begin{equation}\label{e-2.8}
	h_0(z)=\frac{1}{1-z}\exp\left(\frac{2z}{1-z}\right)=\exp\left(\sum_{n=1}^{\infty}\left(2+\frac{1}{n}\right)z^n\right)
\end{equation}
\begin{equation}\label{e-2.9}
	g_0(z)=(1-z)\exp\left(\frac{2z}{1-z}\right)=\exp\left(\sum_{n=1}^{\infty}\left(2-\frac{1}{n}\right)z^n\right).
\end{equation}
Then the function $ f_0 $ defined by 
\begin{equation}\label{e-2.10}
	f_0(z)=zh_0(z)\overline{g_0(z)}=\frac{z(1-\bar{z})}{1-z}\exp\left({\rm Re}\left(\frac{4z}{1-z}\right)\right)\;\; \mbox{for}\;\; z\in\mathbb{D}
\end{equation}
is the log-harmonic Koebe function.\vspace{2mm}

In 2011, Duman \cite{Duman-2011} obtained the upper bound for $ |h(z)| $ and $ |g(z)| $. In $ 2016 $, Ali \textit{et al.} \cite[Theorem 2]{Ali-Abhulhadi-Ng-CVEE-2016} established the sharp lower bounds and exhibited the corresponding extremal functions $ h_0 $, $ g_0 $ and $ f_0 $. Ali \textit{et al.} \cite{Ali-Abhulhadi-Ng-CVEE-2016} extended the Bohr phenomenon to the context of starlike univalent log-harmonic mappings of the form 
\begin{equation}\label{e-22.12}
	f(z)=zh(z)\overline{g(z)}\;\;\text{in}\;\; \mathcal{ST}^{0}_{LH},
\end{equation}
and proved the following interesting result.

\begin{thm}\cite{Ali-Abhulhadi-Ng-CVEE-2016}\label{th-2.3}
	Let $ f $ be a function given by \eqref{e-22.12}. Also, let $ H(z)=zh(z) $ and $ G(z)=zg(z) $. Then
\[
\begin{cases}
	\displaystyle\frac{1}{2e}\leq d(0,\partial H(\mathbb{D}))\leq 1\vspace{3.5mm}\\ \displaystyle\frac{2}{e}\leq d(0,\partial G(\mathbb{D}))\leq 1\vspace{3.5mm}\\	\displaystyle\frac{1}{e^2}\leq d(0,\partial f(\mathbb{D}))\leq 1.
\end{cases}
\]
\vspace{1.5mm}

\noindent Equalities occur if, and only if, $ h $, $ g $ and $ f $ are suitable rotation of $ h_0 $, $ g_0 $ and $ f_0 $. 
\end{thm}

\par In $ 1989 $, Abdulhadi and Hengartner \cite{Abdulhadi-Hengartner-1989} established the sharp coefficient bounds for the function in the class $ \mathcal{ST}^{0}_{LH} $.
\begin{thm}\cite{Abdulhadi-Hengartner-1989}\label{th-2.4}
	Let $ f $ be a function given by \eqref{e-22.12}.  Then 
	\begin{equation*}
		|a_n|\leq 2+\frac{1}{n}\;\; \text{and}\;\; |b_n|\leq 2-\frac{1}{n}\;\; \text{for all}\; n\geq 1.
	\end{equation*}
	Equalities hold for rotation of the function $ f_0 $.
\end{thm}

In $ 2016 $, Ali \emph{et al.} \cite{Ali-Abhulhadi-Ng-CVEE-2016} obtained Bohr radius for log-harmonic mappings of the class $ \mathcal{ST}^{0}_{LH} $.
\begin{thm}\cite{Ali-Abhulhadi-Ng-CVEE-2016}
	Let $ f(z)=zh(z)\overline{g(z)}\in \mathcal{ST}^{0}_{LH}  $ and $ H(z)=zh(z) $ and $ G(z)=zg(z) $. Then
	\begin{enumerate}
		\item[(i)] the inequality 
		\begin{equation*}
			M_h(r):=|z|\exp\left(\sum_{n=1}^{\infty}|a_n||z|^n\right)\leq d(0,\partial H(\mathbb{D}))
		\end{equation*}
		holds for $ |z|\leq r_{H}\approx 0.1222 $, where $ r_H $ is the unique root in $ (0,1) $ of
		\begin{equation*}
			\frac{r}{1-r}\exp\left(\frac{2r}{1-r}\right)=\frac{1}{2e}.
		\end{equation*}
		\item[(ii)] the inequality 
		\begin{equation*}
			M_g(r):=|z|\exp\left(\sum_{n=1}^{\infty}|b_n||z|^n\right)\leq d(0,\partial G(\mathbb{D}))
		\end{equation*}
		holds for $ |z|\leq r_{G}\approx 0.3659 $, where $ r_G $ is the unique root in $ (0,1) $ of
		\begin{equation*}
			r(1-r)\exp\left(\frac{2r}{1-r}\right)=\frac{2}{e}.
		\end{equation*}
	\end{enumerate} 
	Both the radii are sharp and are attained by appropriate rotation of the functions $ H_0(z)=zh_0(z) $  and $ G_0(z)=zg_0(z) $.
\end{thm}

\begin{thm}\cite{Ali-Abhulhadi-Ng-CVEE-2016}
	Let $ f $ be a function given by \eqref{e-22.12}.  Then for any real $ t $, the inequality
	\begin{equation*}
		|z|\exp\left(\sum_{n=1}^{\infty}|a_n+e^{it}b_n||z|^n\right)\leq d(0,\partial f(\mathbb{D}))
	\end{equation*}
	holds for $ |z|\leq r_f\approx0.09078 $, where $ r_f $ is the unique root in $ (0,1) $ of 
	\begin{equation*}
		r\exp\left(\frac{4r}{1-r}\right)=\frac{1}{e^2}.
	\end{equation*}
	The bound is sharp and is attained by suitable rotation of the log-harmonic Koebe function $ f_0 $.
\end{thm}

Our another interest in this paper is to study Bohr radius for the class of analytic functions $ f $ which map unit disk $ \mathbb{D} $ into a concave-wedge domain. The concave-wedge domain is defined (see \cite{Muahnna-ALi-Hasni-2014-JMAA}) by
\begin{equation*}
	W_{\alpha}=\bigg\{w\in\mathbb{C}: |\arg w|<\frac{\alpha\pi}{2},\; 1\leq\alpha\leq 2\bigg\}.
\end{equation*}
It is known that the conformal mapping from $ \mathbb{D} $ onto $ W_{\alpha} $ is given by
\begin{equation}\label{e-2.1111}
	F_{\alpha,t}(z)=t\left(\frac{1+z}{1-z}\right)^{\alpha}=t\left(1+\sum_{n=1}^{\infty}A_nz^n\right) \;\;\mbox{for}\; 1\leq\alpha\leq 2\;\mbox{and}\;\; t>0.
\end{equation}
\par It is easy to see that when $ \alpha=1 $, the domain turns out to be a convex half-plane and when $ \alpha=2 $ it gives a slit domain. Let $ S_{W_{\alpha}} $ be the class of analytic functions $ f $ which maps the unit disk $ \mathbb{D} $ into the wedge domain $ W_{\alpha} $.\vspace{2mm}

\par In $ 2014 $, Abu-Muhana \emph{et al.} \cite{Muahnna-ALi-Hasni-2014-JMAA} proved the following interesting result for functions in the class $ S_{W_{\alpha}} $.
\begin{thm}\label{th-2.13}\cite{Muahnna-ALi-Hasni-2014-JMAA}
	Let $ \alpha\in [1,2] $. If $ f(z)=a_0+\sum_{n=1}^{\infty}a_nz^n\in \mathcal{S}_{W_{\alpha}} $ with $ a_0>0 $, then the inequlaity
	\begin{equation*}
		\sum_{n=1}^{\infty}|a_n||z|^n\leq d(a_0, \partial W_{\alpha})
	\end{equation*}
	holds for $ |z|\leq r_{\alpha}=(2^{1/\alpha}-1)/(2^{1/\alpha}+1) $. The function $ f=F_{\alpha,a_0} $ in \eqref{e-2.1111} shows that $ r_{\alpha} $ is sharp.
\end{thm}
The following lemma is useful to prove one of our main results for functions in the Class $ S_{W_{\alpha}} $.
\begin{lem}\cite{Muahnna-ALi-Hasni-2014-JMAA}
	Let $ F_{\alpha,t} $ be given by \eqref{e-2.1111}, where $ \alpha\in [1,2] $. Then $ A_n>0 $ for all $ n\geq 1 $.
\end{lem}
 \section{Mian results}
 \subsection{Bohr radius in subordination and bounded harmonic classes}
 It is natural to investigate the improved version of the Theorem \ref{th-2.1}. We prove the following improved sharp Bohr radius for the class $ \mathcal{S}(f) $.
\begin{thm}\label{th-2.2}
	Let $\beta \in [0,1/4)$. If $ g(z)= \sum_{n=0}^{\infty}b_nz^n\in \mathcal{S}(f) $ and $ f(z)=\sum_{n=0}^{\infty}a_nz^n $ is univalent, then \begin{equation}\label{e-2.5}
		\beta|f^{\prime}(0)|+\sum_{n=0}^{\infty}|b_nz^n|\leq d(f(0), \partial\Omega)\;\;\mbox{for}\;\; |z|\leq r_{\beta}=\frac{3-4\beta-\sqrt{8}\sqrt{1-2\beta}}{1-4\beta}.
	\end{equation}
The radius $ r_{\beta} $ is sharp for the Koebe function $ f_K(z)=z/(1-z)^2. $
\end{thm}
\begin{rem}
In particular, when $ \beta=0 $, the radius $ r_{\beta} $ which has been proved in Theorem \ref{th-2.2} coincides exactly with $ r_0=3-\sqrt{8}=0.17159 $ in Theorem \ref{th-2.1}. Further, in particular, we obtain $ r_{\beta}=5-2\sqrt{6}\approx 0.10102 $ for $ \beta=1/8 $, $ r_{\beta}=9-4\sqrt{5}\approx 0.05572 $ for $ \beta=3/16  $ and $ r_{\beta}=17-12\sqrt{2}\approx 0.02943 $ for $ \beta=7/32  $. In fact, we see that $ \lim_{\beta\rightarrow ({1}/{4})^{-}}r_{\beta}=0. $
\end{rem}

\subsection{Improved Bohr radius for starlike log-harmonic mappings}
It is known that if $ f $ is a non-vanishing log-harmonic mapping then $ f $ can be written as $ f(z)=h(z)\overline{g(z)} $ where $ h $ and $ g $ are analytic functions in $ \mathbb{D} $. On the other hand, if $ f $ vanishes at $ z=0 $ but is not identically zero, then $ f $ admits the following representation
\begin{equation}\label{e-22.77}
	f(z)=z^m|z|^{2\beta m}h(z)\overline{g(z)}
\end{equation}
where $ m $ is a non-negative integer and $ {\rm Re}\; \beta>-1/2 $, and $ h $, $ g $ are analytic functions in $ \mathbb{D} $ with $ g(0)=1 $ and $ h(0)\neq 1 $ (see \cite{Abdulhadi-Ali-AAA-2012}). The exponent $ \beta $ in \eqref{e-22.77} depends only on $ w(0) $ and it can be expresses as
\begin{equation*}
	\beta=\overline{w(0)}\frac{1+w(0)}{1-|w(0)|^2}.
\end{equation*} 
\par Note that $ f(0)\neq 0 $ if, and only if, $ m=0 $, and that a univalent log-harmonic mapping on $ \mathbb{D} $ vanishes at the origin if, and only if, $ m=1 $. Univalent log-harmonic mappings have been studied extensively by many researchers (see \cite{Ali-Abhulhadi-Ng-CVEE-2016,Duman-2011,Jun-PAMS-1993}). The class of log-harmonic mappings is denoted by $ \mathcal{S}_{LH} $. Let $ z|z|^{2\beta}h(z)\overline{g(z)} $ be a log-harmonic univalent function. We say that $ f $ is a starlike log-harmonic mapping if 
\begin{equation}\label{e-2.7}
	\frac{\partial}{\partial\theta}{\rm Arg}(f(e^{i\theta}))={\rm Re}\left(\frac{zf_z-\bar{z}f_{\bar{z}}}{f}\right)>0\;\; \mbox{in}\;\;\mathbb{D}
\end{equation}
and we denote the set of all strlike log-harmonic functions by $ \mathcal{ST}_{LH} $. Let $ \mathcal{ST}^{0}_{LH} $ be a subclass of $ \mathcal{ST}_{LH} $, consisiting of functions $ f\in\mathcal{ST}_{LH} $ which map the unit disk $ \mathbb{D} $ onto a starlike domain (with repsect to the origin).\vspace{2mm}

\par Our main aim is to study Bohr radius for the class of sense-preserving satrlike log-harmonic mappings in $ \mathbb{D} $ of the form $ f(z)=zh(z)\overline{g(z)}  $ with 
\begin{equation*}
	h(z)=\exp\left(\sum_{n=1}^{\infty}a_kz^k\right)\;\;\mbox{and}\;\; g(z)=\exp\left(\sum_{n=1}^{\infty}b_nz^n\right),
\end{equation*}
where $ h(z) $ and $ g(z) $ may be called as analytic and co-analytic factors of the function $ f(z). $\vspace{2mm}

\par We prove the following improved Bohr radius for functions in the class $ \mathcal{ST}^{0}_{LH}. $
\begin{thm}\label{th-2.7}
		Let $ f $ be a function given by \eqref{e-22.12}. Then for any real $ t $, the inequality
	\begin{equation*}
		|z|\exp\left(\sum_{n=1}^{\infty}\bigg|a_n+e^{it}b_n+\frac{n}{4n^2-1}a_nb_n\bigg||z|^n\right)\leq d(0,\partial f(\mathbb{D}))
	\end{equation*}
	holds for $ |z|\leq r_f\approx0.08528 $, where $ r_f $ is the unique root in $ (0,1) $ of 
	\begin{equation}\label{e-22.13}
		\frac{r}{1-r}\exp\left(\frac{4r}{1-r}\right)=\frac{1}{e^2}\;\;\mbox{in}\;\; (0,1).
	\end{equation}
 The radius $ r_f $ is sharp and is attained by a suitable rotation of the log-harmonic Koebe function $ f_0 $ given by \eqref{e-2.10}.
\end{thm}
\begin{figure}[!htb]
	\begin{center}
		\includegraphics[width=0.50\linewidth]{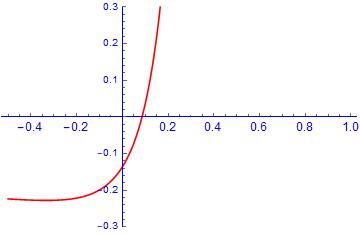}
	\end{center}
	\caption{The radius $ r_f\approx 0.08528$ is a root of \eqref{e-22.13} in $ (0,1) $.}
\end{figure}
\begin{figure}[!htb]
	\begin{center}
		\includegraphics[width=0.46\linewidth]{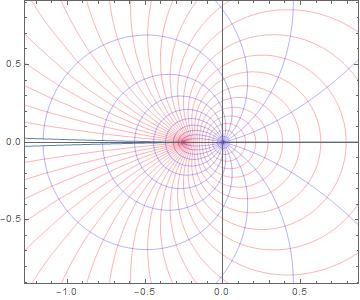}\;\;\;\;\; \includegraphics[width=0.38\linewidth]{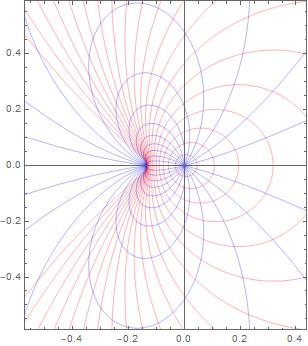}
	\end{center}
	\caption{Image of unit disk $ \mathbb{D} $ under the Koebe function $ f(z)=\frac{z}{(1-z)^2} $ and log-harmonic Koebe function $ f_0(z)=\frac{z(1-\bar{z})}{1-z}\exp\left({\rm Re}\left(\frac{4z}{1-z}\right)\right) $.}
\end{figure}
\begin{figure}[!htb]
	\begin{center}
		\includegraphics[width=0.46\linewidth]{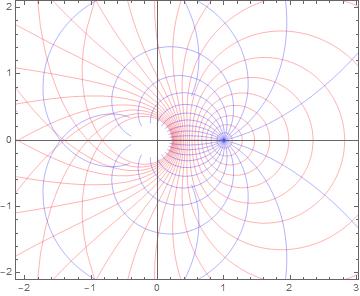}\;\;\;\;\; \includegraphics[width=0.43\linewidth]{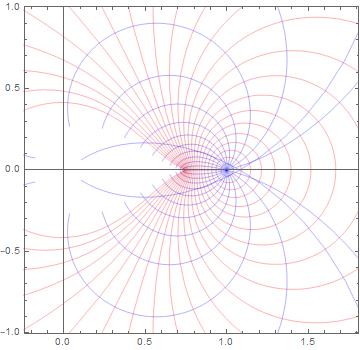}
	\end{center}
	\caption{Image of unit disk $ \mathbb{D} $ under the map $ h_0(z)=\frac{1}{1-z}\exp\left(\frac{2z}{1-z}\right)$ and $ g_0(z)=(1-z)\exp\left(\frac{2z}{1-z}\right)$.}
\end{figure}
\par Next we prove the sharp Bohr radius for the class $ \mathcal{ST}^{0}_{LH} $ in view of additional terms $ |a_n|^2 $ and $ |b_n|^2 $ in the series expansion of $ h $ and $ g $ respectively. 
\begin{thm}\label{th-2.8}
	Let $ f $ be a function given by \eqref{e-22.12} and $ H(z)=zh(z) $ and $ G(z)=zg(z) $. Then
\begin{enumerate}
\item[(i)] the inequality 
\begin{equation*}
|z|\exp\left(\sum_{n=1}^{\infty}\left(|a_n|+\frac{n}{(2n+1)^2}|a_n|^2\right)|z|^n\right)\leq d(0,\partial H(\mathbb{D}))
\end{equation*}
holds for $ |z|\leq r_{H}\approx 0.09735 $, where $ r_H $ is the unique root  of
\begin{equation}\label{e-22.14}
	\frac{r}{(1-r)^2}\exp\left(\frac{2r}{1-r}\right)=\frac{1}{2e} \;\;\mbox{in}\;\; (0,1).
\end{equation}
\begin{figure}[!htb]
	\begin{center}
		\includegraphics[width=0.50\linewidth]{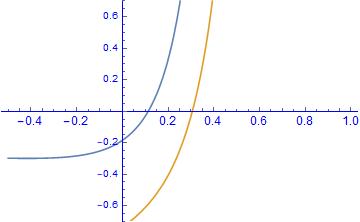}
	\end{center}
	\caption{The radii $ r_H\approx 0.09735$ and $ r_G\approx 0.30539$ are roots of \eqref{e-22.14} and \eqref{e-22.15} respectively in $ (0,1) $.}
\end{figure}
\item[(ii)] the inequality 
\begin{equation*}
|z|\exp\left(\sum_{n=1}^{\infty}\left(|b_n|+\frac{n}{(2n-1)^2}|b_n|^2\right)|z|^n\right)\leq d(0,\partial G(\mathbb{D}))
\end{equation*}
holds for $ |z|\leq r_{G}\approx 0.30539 $, where $ r_G $ is the unique root of
\begin{equation}\label{e-22.15}
r\exp\left(\frac{2r}{1-r}\right)=\frac{2}{e} \;\;\mbox{in}\;\; (0,1).
\end{equation}
\end{enumerate} 
Both the radii are sharp and are attained by appropriate  rotation of $ H_0(z)=zh_0(z) $ and $ G_0(z)=zg_0(z) $.
\end{thm}
We prove the next improved sharp Bohr radius for the class $ \mathcal{ST}^{0}_{LH} $ adding $ |H(z)| $ and $ |G(z)| $ with $ M_h(r) $ and $ M_g(r) $ respecively.
\begin{thm}\label{th-2.12}
	Let $ f $ be a function given by \eqref{e-22.12} and $ H(z)=zh(z) $ and $ G(z)=zg(z) $. Then
	\begin{enumerate}
		\item[(i)]  the inequality 
		\begin{equation*}
			|H(z)|+|z|\exp\left(\sum_{n=1}^{\infty}|a_n||z|^n\right)\leq d(0,\partial H(\mathbb{D}))
		\end{equation*}
		holds for $ |z|\leq r_{H}\approx 0.1073 $, where $ r_H $ is the unqiue root  of
		\begin{equation}\label{e-22.20}
			r\left(\frac{2r}{1-r}-\log(1-r)+\frac{1}{1-r}\exp\left(\frac{2r}{1-r}\right)\right)=\frac{1}{2e} \;\;\mbox{in}\;\; (0,1).
		\end{equation}
		\item[(ii)] the inequality 
		\begin{equation*}
			|G(z)|+|z|\exp\left(\sum_{n=1}^{\infty}|b_n||z|^n\right)\leq d(0,\partial G(\mathbb{D}))
		\end{equation*}
		holds for $ |z|\leq r_{G}\approx 0.3063 $, where $ r_G $ is the unique root of
		\begin{equation}\label{e-22.21}
			r\left(\frac{2r}{1-r}+\log(1-r)+(1-r)\exp\left(\frac{2r}{1-r}\right)\right)=\frac{2}{e} \;\;\mbox{in}\;\; (0,1).
		\end{equation}
	\end{enumerate} 
	Both the radii are sharp and are attained by appropriate rotation of $ H_0(z)=zh_0(z) $ and $ G_0(z)=zg_0(z) $.
\end{thm}
\begin{figure}[!htb]
	\begin{center}
		\includegraphics[width=0.50\linewidth]{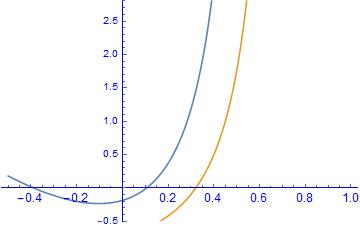}
	\end{center}
	\caption{The radii $ r_H\approx 0.1073$ and $ r_G\approx 0.3063$ are roots of \eqref{e-22.20} and \eqref{e-22.21} respectively in $ (0,1) $.}
\end{figure}
\par For any positive integer $ m $, considering $ |h(z)|^m $ and $ |g(z)|^m $, next we prove the improved sharp Bohr radius for the class $ \mathcal{ST}^{0}_{LH} $. 
\begin{thm}\label{th-2.9}
	Let $ f $ be a function given by \eqref{e-22.12} and $ H(z)=zh(z) $ and $ G(z)=zg(z) $. 
	\begin{enumerate}
		\item[(i)] If $ |h(z)|\leq 1 $, then for any $ m\in\mathbb{N} $, the inequality 
		\begin{equation*}
			|z|\exp\left(|h(z)|^m+\sum_{n=1}^{\infty}|a_n||z|^n\right)\leq d(0,\partial H(\mathbb{D}))
		\end{equation*}
		holds for $ |z|\leq r_{H}\approx 0.0566 $, where $ r_H $ is the unique root  of
		\begin{equation}\label{e-22.16}
			\frac{re}{1-r}\exp\left(\frac{2r}{1-r}\right)=\frac{1}{2e} \;\;\mbox{in}\;\; (0,1).
		\end{equation}
		\item[(ii)] If $ |g(z)|\leq 1 $, then for any $ m\in\mathbb{N} $, the inequality 
		\begin{equation*}
			|z|\exp\left(|g(z)|^m+\sum_{n=1}^{\infty}|b_n||z|^n\right)\leq d(0,\partial G(\mathbb{D}))
		\end{equation*}
		holds for $ |z|\leq r_{G}\approx 0.1764 $, where $ r_G $ is the unique root of
		\begin{equation}\label{e-22.17}
			re(1-r)\exp\left(\frac{2r}{1-r}\right)=\frac{2}{e} \;\;\mbox{in}\;\; (0,1).
		\end{equation}
	\end{enumerate} 
	Both the radii are sharp and are attained by a suitable rotation of $ H_0(z)=zh_0(z) $ and $ G_0(z)=zg_0(z) $.
	\begin{figure}[!htb]
		\begin{center}
			\includegraphics[width=0.50\linewidth]{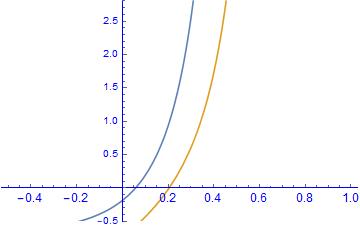}
		\end{center}
		\caption{The radii $ r_H\approx 0.0566$ and $ r_G\approx 0.1764$ are roots of \eqref{e-22.16} and \eqref{e-22.17}  respectively in $ (0,1) $.}
	\end{figure}
\end{thm}
\begin{rem}
	It is worth to notice in Theorem \ref{th-2.9} that the Bohr radii $ r_H $ and $ r_G $ are independent of the choice of the positive integer $ m $.
\end{rem}
We prove the improved sharp Bohr radius adding $ |h(z)+g(z)| $ with the series  $ \sum_{n=1}^{\infty}|a_n+e^{it}b_n||z|^n $ for the class $ \mathcal{ST}^{0}_{LH} $.
\begin{thm}\label{th-2.10}
Let $ f $ be a function given by \eqref{e-22.12} with $ |h(z)|+|g(z)|\leq 1 $. Then for any real $ t $, the inequality
\begin{equation*}
|z|\exp\left(|h(z)+g(z)|+\sum_{n=1}^{\infty}|a_n+e^{it}b_n||z|^n\right)\leq d(0,\partial f(\mathbb{D}))
\end{equation*}
holds for $ |z|\leq r_f\approx0.04181 $, where $ r_f $ is the unique root of 
\begin{equation}\label{e-22.18}
er\exp\left(\frac{4r}{1-r}\right)=\frac{1}{e^2} \;\;\mbox{in}\;\; (0,1).
\end{equation}
The Bohr radius $ r_f $ is sharp and is attained by suitable rotation of the log-harmonic Koebe function $ f_0 $.
\end{thm}
Next we prove the improved sharp Bohr radius for the class $ \mathcal{ST}^{0}_{LH} $ adding $ |f(z)| $.
\begin{thm}\label{th-2.11}
Let $ f $ be a function given by \eqref{e-22.12} with $ |h(z)|\leq 1 $ and $|g(z)|\leq 1 $. Then for any real $ t $, the inequality
\begin{equation*}
|f(z)|+|z|\exp\left(\sum_{n=1}^{\infty}|a_n+e^{it}b_n||z|^n\right)\leq d(0,\partial f(\mathbb{D}))
\end{equation*}
holds for $ |z|\leq r_f\approx0.0592 $, where $ r_f $ is the unique root of 
\begin{equation}\label{e-22.19}
r\left(1+\exp\left(\frac{4r}{1-r}\right)\right)=\frac{1}{e^2} \;\;\mbox{in}\;\; (0,1).
\end{equation}
The Bohr radius $ r_f $ is sharp for a suitable rotation of the log-harmonic Koebe function $ f_0 $.
\end{thm}
\begin{figure}[!htb]
\begin{center}
\includegraphics[width=0.50\linewidth]{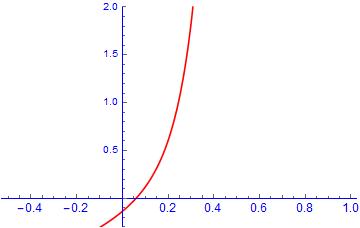}
\end{center}
\caption{The radius $ r_f\approx 0.0592$ is the root of \eqref{e-22.19} in $ (0,1) $.}
\end{figure}
\subsection{Bohr radius for concave-wedge domain}
We prove the following result which is an improved version of Theorem \ref{th-2.13}.
\begin{thm}\label{th-2.14}
Let $ \alpha\in [1,2] $ and $ \beta\in [0,2) $. If $ f(z)=a_0+\sum_{n=1}^{\infty}a_nz^n\in S_{W_{\alpha}} $ with $ a_0>0 $, then the inequlaity
\begin{equation*}
\frac{2\beta a_0}{\alpha\pi}|\arg f(z)|+\sum_{n=1}^{\infty}|a_n||z|^n\leq d(a_0, \partial W_{\alpha})
\end{equation*}
holds for $ |z|\leq r_{\alpha,\beta}=((2-\beta)^{1/\alpha}-1)/((2-\beta)^{1/\alpha}+1) $. The function $ f=F_{\alpha,a_0} $ in \eqref{e-2.1111} shows that $ r_{\alpha,\beta} $ is sharp.
\end{thm}
\begin{rem}
Since $ W_{\alpha} $ turns out to be a convex half-plane when $ \alpha=1 $, it is evident that, for $ \alpha=1 $ and $ \beta=0 $, the radius $ r_{\alpha, \beta} $ coincides exactly with the Bohr radius $ 1/3 $.
\end{rem}
\begin{figure}[!htb]
\begin{center}
\includegraphics[width=0.25\linewidth]{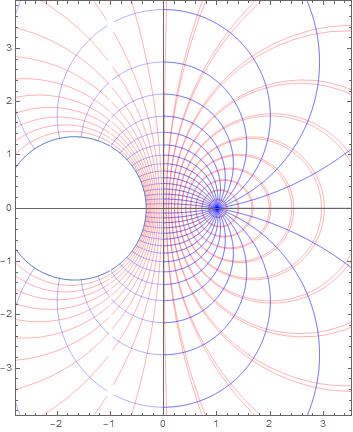}\;\;\; \includegraphics[width=0.364\linewidth]{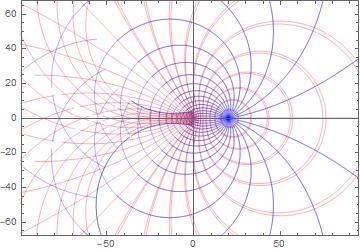}\;\;\;\includegraphics[width=0.30\linewidth]{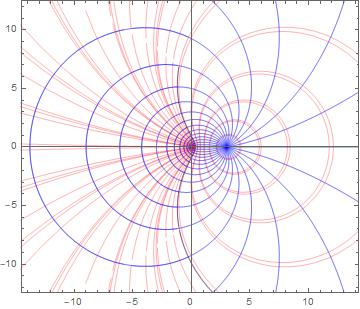}
\end{center}
\caption{Image of unit disk $ \mathbb{D} $ under the maps $ F_{1,1}(z) $, $ F_{1.5,20}(z) $ and $ F_{2,3}(z) $ repsectively.}
\end{figure}
\section{Proof of the main results}
\begin{proof}[\bf Proof of Theorem \ref{th-2.2}]
	Since $ g(z)=\sum_{n=0}^{\infty}b_nz^n\in S(f) $ and $ f(z)=\sum_{n=0}^{\infty}a_nz^n $ are univalent, by the famous well-known de Brandge's theorem \cite[p.17]{Gong-1991}, we have
	\begin{equation}\label{e-3.1}
		|b_n|\leq n|f^{\prime}(0)|.
	\end{equation}
Therefore, from \eqref{e-2.4} and \eqref{e-3.1}, it is easy to see that
\begin{equation}\label{e-3.2}
	|b_n|\leq 4nd(f(0),\partial\Omega).
\end{equation}
By a simple computation using \eqref{e-3.2}, we obtain
\begin{align*}
	\beta|f^{\prime}(0)|+\sum_{n=1}^{\infty}|b_nz^n|&=\beta|f^{\prime}(0)|+\sum_{n=1}^{\infty}|b_n||z|^n\\&\leq 4\beta d(f(0),\partial\Omega)+4d(f(0),\partial\Omega)\sum_{n=1}^{\infty}nr^n\\&=4d(f(0),\partial\Omega)\left(\beta+\frac{r}{(1-r)^2}\right)\\&\leq d(f(0),\partial\Omega)
\end{align*}
if, and only if,
\begin{equation*}
	4\left(\beta+\frac{r}{(1-r)^2}\right)\leq 1.
\end{equation*}
Therefore, \eqref{e-2.5} holds for 
\begin{equation*}
	|z|\leq r_0=\frac{3-4\beta-\sqrt{8}\sqrt{1-2\beta}}{1-4\beta}.
\end{equation*}
Since 
\begin{equation*}
	f_K(z)=\frac{z}{(1-z)^2}=\sum_{n=1}^{\infty}nz^n\;\; \text{and}\;\; d(f(0). \partial\Omega)=\frac{1}{4},
\end{equation*}
 a simple computation shows that
 \begin{align*}
 	\beta|f_K^{\prime}(0)|+\sum_{n=1}^{\infty}|b_nz^n|&=\beta|f_K^{\prime}(0)|+\sum_{n=1}^{\infty}|b_n|r_0^n\\&= 4\beta d(f_K(0),\partial\Omega)+4d(f_K(0),\partial\Omega)\sum_{n=1}^{\infty}nr_0^n\\&=4d(f_K(0),\partial\Omega)\left(\beta+\frac{r_o}{(1-r_0)^2}\right)\\&= d(f(0),\partial\Omega).
 \end{align*}
This shows that the radius $ r_0 $ is the best possible. This completes the proof.
\end{proof}	
\begin{proof}[\bf Proof of Theorem \ref{th-2.7}]
Let $ f(z)=zh(z)\overline{g(z)}\in \mathcal{ST}^{0}_{LH}  $. Then in view of Theorem \ref{th-2.4}, we have the following coefficient bounds
\begin{equation*}
	|a_n|\leq 2+\frac{1}{n}\;\; \text{and}\;\; |b_n|\leq 2-\frac{1}{n}\;\; \text{for all}\; n\geq 1.
\end{equation*}
On the other hand, from Theorem \ref{th-2.3}, we have
\begin{equation*}
	d(0,\partial f(\mathbb{D}))\geq \frac{1}{e^2}.
\end{equation*}
Therefore, a simple computation shows that
\begin{align*} &
	r\exp\left(\sum_{n=1}^{\infty}|a_n|r^n+\sum_{n=1}^{\infty}|b_n|r^n+\sum_{n=1}^{\infty}\frac{n}{4n^2-1}|a_n||b_n|r^n\right)\\&\leq r \exp\left(\sum_{n=1}^{\infty}\left(2+\frac{1}{n}\right)r^n+\sum_{n=1}^{\infty}\left(2-\frac{1}{n}\right)r^n+\sum_{n=1}^{\infty}\frac{n}{4n^2-1}\left(2+\frac{1}{n}\right)\left(2-\frac{1}{n}\right)r^n\right)\\&= r\exp\left(4\sum_{n=1}^{\infty}r^n+\sum_{n=1}^{\infty}\frac{r^n}{n}\right)\\&=r\exp\left(\frac{4r}{1-r}-\log(1-r)\right)\\&\leq d(0,\partial f(\mathbb{D}))
\end{align*}
if, and only if,
\begin{equation*}
	r\exp\left(\frac{4r}{1-r}-\log(1-r)\right)\leq\frac{1}{e^2}
\end{equation*}
which is equivalent to
\begin{equation*}
	\frac{r}{1-r}\exp\left(\frac{4r}{1-r}\right)\leq\frac{1}{e^2}.
\end{equation*}
The Bohr radius $ r_f $ is the unique root of the equation
\begin{equation*}
	\frac{r}{1-r}\exp\left(\frac{4r}{1-r}\right)=\frac{1}{e^2}
\end{equation*}
in $ (0,1), $ a simple com putation shows that $ r_f\approx 0.08528 $.\vspace{1mm}

\noindent In order to show the sharpness of $ r_f $, let $ h_0 $, $ g_0 $ and $ f_0 $ be given by \eqref{e-2.8}, \eqref{e-2.9} and \eqref{e-2.10} repsectively. For these functions, it is easy to see that 
\begin{equation}\label{e-3.3}
	|a_n|=2+\frac{1}{n},\;\; |b_n|=2-\frac{1}{n}\;\; \text{for all}\; n\in\mathbb{N}\;\;\text{and}\;\; d(0,\partial f_0(\mathbb{D}))=\frac{1}{e^2}.
\end{equation}
A simple computation using \eqref{e-3.3} shows that
\begin{align*} &
	r_f\exp\left(\sum_{n=1}^{\infty}|a_n|r_f^n+\sum_{n=1}^{\infty}|b_n|r_f^n+\sum_{n=1}^{\infty}\frac{n}{4n^2-1}|a_n||b_n|r_f^n\right)\\&= r_f \exp\left(\sum_{n=1}^{\infty}\left(2+\frac{1}{n}\right)r_f^n+\sum_{n=1}^{\infty}\left(2-\frac{1}{n}\right)r_f^n+\sum_{n=1}^{\infty}\frac{n}{4n^2-1}\left(2+\frac{1}{n}\right)\left(2-\frac{1}{n}\right)r_f^n\right)\\&= r_f\exp\left(4\sum_{n=1}^{\infty}r_f^n+\sum_{n=1}^{\infty}\frac{r_f^n}{n}\right)\\&=r_f\exp\left(\frac{4r_f}{1-r_f}-\log(1-r_f)\right)\\&=\frac{1}{e^2}\\&= d(0,\partial f_0(\mathbb{D})).
\end{align*}
Therefore, the radius $ r_f $ is the best possible.
\end{proof}
\begin{proof}[\bf Proof of Theorem \ref{th-2.8}]
	Let $ f(z)=zh(z)\overline{g(z)}\in \mathcal{ST}^{0}_{LH}  $, $ H(z)=zh(z) $ and $ G(z)=zg(z). $ (i) In view of  Theorem \ref{th-2.3} and Theorem \ref{th-2.4}, we have the following sharp coefficient estimates
	\begin{equation*}
		|a_n|\leq 2+\frac{1}{n}\;\; \text{for all}\;\; n\geq 1\;\; \text{and}\;\; d(0,\partial H(\mathbb{D}))\geq\frac{1}{2e}
	\end{equation*}
and equality holds for the function $ h_0(z) $ and $ H(z)=zh_0(z) $ respectively. A simple computation shows that
\begin{align*}&
	|z|\exp\left(\sum_{n=1}^{\infty}\left(|a_n|+\frac{n}{(2n+1)^2}|a_n|^2\right)|z|^n\right)\\&\leq r\exp\left(2\sum_{n=1}^{\infty}r^n+2\sum_{n=1}^{\infty}\frac{r^n}{n}\right)\\&= r\exp\left(\frac{2r}{1-r}-2\log(1-r)\right)\\&=\frac{r}{(1-r)^2}\exp\left(\frac{2r}{1-r}\right)\\&\leq d(0,\partial H(\mathbb{D}))
\end{align*}
if, and only if,
\begin{equation*}
	\frac{r}{(1-r)^2}\exp\left(\frac{2r}{1-r}\right)\leq\frac{1}{2e}.
\end{equation*}
Therefore, the Bohr radius $ r_H $ is the unique root of the equation
\begin{equation*}
	\frac{r}{(1-r)^2}\exp\left(\frac{2r}{1-r}\right)=\frac{1}{2e}
\end{equation*}
in $ (0,1) $, a computation shows that $ r_H\approx 0.09735. $\vspace{1mm}

\noindent To show the sharpness of the radius $ r_H $, let $ H_0(z)=zh_0(z) $, where $ h_0 $ is given by \eqref{e-2.8}. It is easy to see that 
\begin{align*} &
	|z|\exp\left(\sum_{n=1}^{\infty}\left(|a_n|+\frac{n}{(2n+1)^2}|a_n|^2\right)|z|^n\right)\\&= r_H\exp\left(2\sum_{n=1}^{\infty}r_H^n+2\sum_{n=1}^{\infty}\frac{r_H^n}{n}\right)\\&= r\exp\left(\frac{2r}{1-r_H}-2\log(1-r_H)\right)\\\vspace{2mm}&=	\frac{r_H}{(1-r_H)^2}\exp\left(\frac{2r_H}{1-r_H}\right)\\&\vspace{2mm}=\frac{1}{2e}\\\vspace{2mm}&\vspace{2mm}= d(0,\partial H_0(\mathbb{D})).
\end{align*}
Therefore, the radius $ r_H $ is the best possible.\vspace{1mm}

\noindent (ii) In view of Theorem \ref{th-2.3} and Theorem \ref{th-2.4}, we have
\begin{equation*}
	|b_n|\leq 2-\frac{1}{n}\;\; \text{for all}\;\; n\geq 1\;\; \text{and}\;\; d(0,\partial G(\mathbb{D}))\geq\frac{2}{e}.
\end{equation*}
Both the equalities hold  for the function $ g_0(z) $ and $ G_0(z)=zg_0(z) $ respectively.\vspace{1mm} 

\noindent A simple computation shows that
\begin{align*} &
	|z|\exp\left(\sum_{n=1}^{\infty}\left(|b_n|+\frac{n}{(2n-1)^2}|b_n|^2\right)|z|^n\right)\\&\leq r\exp\left(2\sum_{n=1}^{\infty}r^n\right)\\&= r\exp\left(\frac{2r}{1-r}\right)\\&\leq d(0,\partial G(\mathbb{D}))
\end{align*}
if, and only if,
\begin{equation*}
	r\exp\left(\frac{2r}{1-r}\right)\leq\frac{2}{e}.
\end{equation*}
Therefore, the Bohr radius $ r_G $ is the unique root of the equation
\begin{equation*}
		r\exp\left(\frac{2r}{1-r}\right)=\frac{2}{e}
\end{equation*}
in $ (0,1) $ which can be computed as $ r_G\approx 0.30539. $\vspace{1mm}

\noindent To show the sharpness of the radius $ r_G $, let $ G_0(z)=zg_0(z) $, where $ g_0 $ is given by \eqref{e-2.8}. It is easy to see that 
\begin{align*} &
	|z|\exp\left(\sum_{n=1}^{\infty}\left(|b_n|+\frac{n}{(2n-1)^2}|b_n|^2\right)|z|^n\right)\\&= r_G\exp\left(2\sum_{n=1}^{\infty}r_G^n\right)\\&= r_G\exp\left(\frac{2r_G}{1-r_G}\right)\\&=\frac{2}{e}\\&= d(0,\partial G_0(\mathbb{D})).
\end{align*}
Therefore, the radius $ r_G $ is the best possible.
\end{proof}
\begin{proof}[\bf Proof of Theorem \ref{th-2.9}]
Let $ f(z)=zh(z)\overline{g(z)}\in \mathcal{ST}^{0}_{LH}  $, $ H(z)=zh(z) $ and $ G(z)=zg(z) $.\vspace{1mm}

\noindent (i) In view of Theorem \ref{th-2.3} and Theorem \ref{th-2.4} and using the fact that $ |h(z)|\leq 1 $, for $ m\in\mathbb{N} $, we obtain
\begin{align*}
	|z|\exp\left(|h(z)|^m+\sum_{n=1}^{\infty}|a_n||z|^n\right)&= |z|\exp(|h(z)|^m)\exp\left(\sum_{n=1}^{\infty}|a_n||z|^n\right)\\&\leq re \exp\left(\frac{2r}{1-r}-2\log (1-r)\right)\\&=	\frac{re}{1-r}\exp\left(\frac{2r}{1-r}\right)\\&\leq d(0,\partial H(\mathbb{D}))
\end{align*}
if, and only if,
\begin{equation*}
	\frac{re}{1-r}\exp\left(\frac{2r}{1-r}\right)\leq \frac{1}{2e}.
\end{equation*}
Therefore, the Bohr radius $ r_H $ is the solution of 
\begin{equation*}
	\frac{er}{1-r}\exp\left(\frac{2r}{1-r}\right)= \frac{1}{2e}.
\end{equation*}
It is easy to see that $ r_H\approx 0.0566 $. The radius $ r_H $ is best possible and it can be shown by using the function $ H_0(z)=zh_0(z) $. \vspace{1mm}

\noindent (ii) Since $ |g(z)|\leq 1 $, in view of Theorem \ref{th-2.3} and Theorem \ref{th-2.4}, by a simple computation, we obtain
\begin{align*}
	 |z|\exp\left(|g(z)|^m+\sum_{n=1}^{\infty}|b_n||z|^n\right)=& |z|\exp\left(|g(z)|^m\right)\exp\left(\sum_{n=1}^{\infty}|b_n||z|^n\right)\\&\leq er\exp\left(\sum_{n=1}^{\infty}\left(2-\frac{1}{n}\right)r^n\right)\\&=re \exp\left(\frac{2r}{1-r}+\log(1-r)\right)\\&=re(1-r)\exp\left(\frac{2r}{1-r}\right)\\&\leq d(0,\partial G(\mathbb{D}))
\end{align*}
if, and only if,
\begin{equation*}
	re(1-r)\exp\left(\frac{2r}{1-r}\right)\leq\frac{2}{e}.
\end{equation*}
Therefore, the Bohr radius $ r_G $ is the unieque root of the equation 
\begin{equation*}
	re(1-r)\exp\left(\frac{2r}{1-r}\right)=\frac{2}{e}.
\end{equation*}
\noindent  A simple computation shows $ r_G\approx 0.1764 $. The radius $ r_G $ is the best possible which can be shown by considering the function $ G_0(z)=zg_0(z) $, where $ g_0(z) $ is defined in \eqref{e-2.9}.
\end{proof}
\begin{proof}[\bf Proof of Theorem \ref{th-2.10}]
In view of Theorem \ref{th-2.4}, we have the following sharp coefficient bounds
	\begin{equation*}
		|a_n|\leq 2+\frac{1}{n}\;\; \text{and}\;\; |b_n|\leq 2-\frac{1}{n}\;\; \text{for all}\; n\geq 1,
	\end{equation*}
which are attained by the function $ h_0 $ and $ g_0 $ respectively defined in \eqref{e-2.8} and \eqref{e-2.9}. On the other hand, by Theorem \ref{th-2.3}, we have the sharp distance
\begin{equation*}
	d(0,\partial f(\mathbb{D}))\geq \frac{1}{e^2}
\end{equation*}
which is attained by the function $ f_0 $ defined in \eqref{e-2.10}.

\noindent Since $ |h(z)|+|g(z)|\leq 1 $, we obtain
\begin{align*}
	&|z|\exp\left(|h(z)+g(z)|+\sum_{n=1}^{\infty}|a_n||z|^n+\sum_{n=1}^{\infty}|b_n||z|^n\right)\\&= r\exp\left(|h(z)|+|g(z)|\right)\exp\left(\sum_{n=1}^{\infty}|a_n|r^n+\sum_{n=1}^{\infty}|b_n|r^n\right)\\&\leq er \exp\left(4\sum_{n=1}^{\infty}|a_n|r^n\right)\\&=er\exp\left(\frac{4r}{1-r}\right)\\&\leq d(0,\partial f(\mathbb{D}))
\end{align*}
if, and only if,
\begin{equation*}
	er\exp\left(\frac{4r}{1-r}\right)\leq \frac{1}{e^2}.
\end{equation*}
Therefore, the Bohr radius $ r_f $ is the unique root of the equation
\begin{equation*}
	er\exp\left(\frac{4r}{1-r}\right)=\frac{1}{e^2}
\end{equation*}
in $ (0,1) $ which yields $ r_f\approx 0.04181 $. The sharpness of the radius $ r_f $ can be shown by considering a suitable rotation of the log-harmonic Koebe function $ f_0 $. This completes the proof.
\end{proof}
\begin{proof}[\bf Proof of Theorem \ref{th-2.11}]
	Since $ |h(z)|\leq 1 $ and $ |g(z)|\leq 1 $, in view of Theorem \ref{th-2.3} and Theorem \ref{th-2.4}, we obtain
	\begin{align*}&
		|f(z)|+|z|\exp\left(\sum_{n=1}^{\infty}|a_n||z|^n+\sum_{n=1}^{\infty}|b_n||z|^n\right)\\&\leq r|h(z)||g(z)|+r\exp\left(\sum_{n=1}^{\infty}|a_n|r^n+\sum_{n=1}^{\infty}|b_n|r^n\right)\\&\leq r+r\exp\left(4\sum_{n=1}^{\infty}\frac{r^n}{r}\right)\\&=r+r\exp\left(\frac{4r}{1-r}\right)\\&\leq d(0,\partial f(\mathbb{D}))
	\end{align*}
if, and only if,
\begin{equation*}
	r\left(1+\exp\left(\frac{4r}{1-r}\right)\right)\leq \frac{1}{e^2}.
\end{equation*}
Thus, the Bohr radius $ r_f $ is the unique root of the equation 
\begin{equation*}
	r\left(1+\exp\left(\frac{4r}{1-r}\right)\right)= \frac{1}{e^2}
\end{equation*}
which gives $ r_f\approx 0.0592 $. The radius $ r_f $ is sharp and can be shown by considering a suitable rotation of the log-harmonic Koebe function $ f_0. $ This completes the proof.
\end{proof}
\begin{proof}[\bf Proof of Theorem \ref{th-2.12}]
	By Theorem \ref{th-2.3}, we have 
	\begin{equation}\label{e-3.4}
		d(0,\partial H(\mathbb{D}))\geq\frac{1}{2e}\;\;\text{and}\;\; 	d(0,\partial G(\mathbb{D}))\geq\frac{2}{e}
	\end{equation}
and the sharp coefficient bounds
\begin{equation*}
	|a_n|\leq 2+\frac{1}{n}\;\; \text{and}\;\; |b_n|\leq 2-\frac{1}{n}\;\; \text{for all}\; n\geq 1.
\end{equation*}
All the inequlities are attained by the extremal functions $ H_0(z)=zh_0(z) $ and $ G_0(z)=zg_0(z) $, where $ h_0 $ and $ g_0 $ are defined respectively in \eqref{e-2.8} and \eqref{e-2.9}.\vspace{1mm}

\noindent (i) Using Theorem \ref{th-2.4} and \eqref{e-3.4}, we obtain
\begin{align*}
	|H(z)|+|z|\exp\left(\sum_{n=1}^{\infty}|a_n||z|^n\right)&\leq |z|\left(\sum_{n=1}^{\infty}|a_n||z|^n\right)+|z|\exp\left(\sum_{n=1}^{\infty}|a_n||z|^n\right)\\&\leq r\left(\sum_{n=1}^{\infty}\left(2+\frac{1}{n}\right)r^n\right)+r\exp\left(\sum_{n=1}^{\infty}\left(2+\frac{1}{n}\right)r^n\right)\\&=r\left(\frac{2r}{1-r}-\log(1-r)\right)+r\exp\left(\frac{2r}{1-r}-\log(1-r)\right)\\&\leq d(0,\partial H(\mathbb{D}))
\end{align*}
if, and only if,
\begin{equation*}
	r\left(\frac{2r}{1-r}-\log(1-r)\right)+r\exp\left(\frac{2r}{1-r}-\log(1-r)\right)\leq\frac{1}{2e}.
\end{equation*}
Therefore, the Bohr radius $ r_H $ is the unique root of the equation \begin{equation*}
	r\left(\frac{2r}{1-r}-\log(1-r)\right)+r\exp\left(\frac{2r}{1-r}-\log(1-r)\right)=\frac{1}{2e}
\end{equation*}
in $ (0,1) $ which shows that $ r_H\approx 0.1073 $. In order to show the radius $ r_H $ is sharp, we consider the function $ H_0(z)=zh_0(z) $, where $ h_0 $ is defined in \eqref{e-2.8}. Therefore, a simple compuatation shows that
\begin{align*}&
	|H_0(z)|+|z|\exp\left(\sum_{n=1}^{\infty}|a_n||z|^n\right)\\&= |z|\left(\sum_{n=1}^{\infty}|a_n||z|^n\right)+|z|\exp\left(\sum_{n=1}^{\infty}|a_n||z|^n\right)\\&= r_H\left(\sum_{n=1}^{\infty}\left(2+\frac{1}{n}\right)r_H^n\right)+r_H\exp\left(\sum_{n=1}^{\infty}\left(2+\frac{1}{n}\right)r_H^n\right)\\&=r_H\left(\frac{2r_H}{1-r_H}-\log(1-r_H)\right)+r_H\exp\left(\frac{2r_H}{1-r_H}-\log(1-r_H)\right)\\&=\frac{1}{2e}\vspace{2mm}\\&=d(0,\partial H(\mathbb{D})).
\end{align*}
This shows that $ r_H $ is the best possible.\vspace{1.5mm}

\noindent (ii) Using Theorem \ref{th-2.4} and \eqref{e-3.4}, we obtain
\begin{align*}&
	|G(z)|+|z|\exp\left(\sum_{n=1}^{\infty}|b_n||z|^n\right)\\&\leq |z|\left(\sum_{n=1}^{\infty}|b_n||z|^n\right)+|z|\exp\left(\sum_{n=1}^{\infty}|b_n||z|^n\right)\\&\leq r\left(\sum_{n=1}^{\infty}\left(2-\frac{1}{n}\right)r^n\right)+r\exp\left(\sum_{n=1}^{\infty}\left(2+\frac{1}{n}\right)r^n\right)\\&=r\left(\frac{2r}{1-r}+\log(1-r)\right)+r\exp\left(\frac{2r}{1-r}+\log(1-r)\right)\\&\leq d(0,\partial H(\mathbb{D}))
\end{align*}
if, and only if,
\begin{equation*}
	r\left(\frac{2r}{1-r}+\log(1-r)\right)+r\exp\left(\frac{2r}{1-r}+\log(1-r)\right)\leq\frac{2}{e}.
\end{equation*}
Therefore, the Bohr radius $ r_G $ is the unique root of the equation \begin{equation*}
	r\left(\frac{2r}{1-r}-\log(1-r)\right)+r\exp\left(\frac{2r}{1-r}-\log(1-r)\right)=\frac{2}{e}
\end{equation*}
in $ (0,1) $ which yeilds $ r_G\approx 0.3063 $. To show the radius $ r_G $ is best possible, we consider the function $ G_0(z)=zg_0(z) $, where $ g_0 $ is defined in \eqref{e-2.9}.\vspace{2mm}

\noindent Thus, it is easy to see that
\begin{align*}&
	|G_0(z)|+|z|\exp\left(\sum_{n=1}^{\infty}|b_n||z|^n\right)\\&= |z|\left(\sum_{n=1}^{\infty}|b_n||z|^n\right)+|z|\exp\left(\sum_{n=1}^{\infty}|b_n||z|^n\right)\\&= r_G\left(\sum_{n=1}^{\infty}\left(2-\frac{1}{n}\right)r_G^n\right)+r_G\exp\left(\sum_{n=1}^{\infty}\left(2+\frac{1}{n}\right)r_G^n\right)\\&=r_G\left(\frac{2r_G}{1-r_G}+\log(1-r_G)\right)+r_G\exp\left(\frac{2r_G}{1-r_G}+\log(1-r_G)\right)\\&=\frac{2}{e}\\\vspace{2mm}&= d(0,\partial G(\mathbb{D}))
\end{align*}
which shows that $ r_G $ is sharp.
\end{proof}

\begin{proof}[\bf Proof of Theorem \ref{th-2.14}]
	Since $ f\in S_{W_{\alpha}} $, we have $ |\arg f(z)|\leq \pi\alpha/2 $. Therefore, in view of Lemma \ref{lem-2.3}, it is easy to see that
	\begin{align*}
		\frac{2\beta a_0}{\pi\alpha}|\arg f(z)|+\sum_{n=1}^{\infty}|a_nz^n|&\leq \frac{2\beta a_0}{\pi\alpha}\frac{\pi\alpha}{2}+a_0\sum_{n=1}^{\infty}A_nr^n\\&=a_0\left(\sum_{n=1}^{\infty}A_nr^n+\beta\right)\\&=a_0\left(\left(\frac{1+r}{1-r}\right)^{\alpha}-1+\beta\right)\\&= d(a_0,\partial W_{\alpha})\left(\left(\frac{1+r}{1-r}\right)^{\alpha}-1+\beta\right)\\&\leq d(a_0,\partial W_{\alpha})
	\end{align*}
if, and only if,
\begin{equation*}
	\left(\frac{1+r}{1-r}\right)^{\alpha}-1+\beta\leq 1.
\end{equation*}
Therefore, the Bohr radius is the unique root of the equation 
\begin{equation*}
	\left(\frac{1+r}{1-r}\right)^{\alpha}-1+\beta=1.
\end{equation*}
A simple computation shows that $ r_{\alpha,\beta}=((2-\beta)^{1/\alpha}-1)/((2-\beta)^{1/\alpha}+1) $. By a suitable rotation of the function $ f=F_{\alpha,a_0} $ in \eqref{e-2.1111}, it can be shown that the radius $ r_{\alpha,\beta} $ is sharp. This completes the proof.
\end{proof}

\noindent\textbf{Acknowledgment:}  The first author is supported by the Institute Post Doctoral Fellowship of IIT Bhubaneswar, India, the second author is supported by SERB-MATRICS, India.





\begin{thebibliography}{99}	
	
\bibitem{Abdulhadi-Ali-AAA-2012} {\sc Z. Abdulhadi} and {\sc R. M. Ali}, Univalent logharmonic mappings in the plane, \textit{Abstr. Appl. Anal.} 2012: 2012: 32.


\bibitem{Abdulhadi-Hengartner-1989} {\sc Z. Abdulhadi} and {\sc W. Hengartner}, Univalent logharmonic mappings on the left half-plane with periodic dilations, In: H. M. Srivastava, S. Owa, editors. Univalent functions, fractional
calculus, and their applications. Ellis Horwood series in Mathematics and Applications, Chick-
ester: Horwood: 1989, 13-28.

\bibitem{Muhanna-CVEE-2010} {\sc Y. Abu-Muhanna}, Bohr's phenomenon in subordination and bounded harmonic classes,
\textit{Complex Var. Elliptic Equ.} \textbf{55} (2010), 1071-1078.

\bibitem{Muahnna-ALi-Hasni-2014-JMAA} {\sc Y. Abu-Muhanna, R.M. Ali} and {\sc S. M. Hasni}, Bohr radius for subordinating families of analytic functions andbounded harmonic mappings, \textit{J. Math. Anal. Appl.} \textbf{420}(2014), 124-136.

\bibitem{Aizenberg-PAMS-2000} {\sc L. Aizenberg}, Multidimensional analogues of Bohr's theorem on power series, \textit{Proc. Amer. Math. Soc.} \textbf{128} (2000), 1147-1155.

\bibitem{Aizenberg-Tarkhanov-PLMS-2001} {\sc L. Aizenberg} and {\sc N. Tarkhanov}, A Bohr phenomenon for elliptic equations, \textit{Proc. London Math. Soc.} \textbf{82}(2) (2001), 385-401.

\bibitem{Aizenberg-Aytuna-Djakov-PAMS-2000} {\sc L. Aizenberg, A. Aytuna} and {\sc P. Djakov}, An abstract approach to Bohr's Theorem, \textit{Proc. Amer. Math. Soc.} \textbf{128}(9) (2000), 2611-2619.

\bibitem{Aizenberg-Aytuna-JMAA-2001} {\sc L. Aizenberg, A. Aytuna} and {\sc P. Djakov}, Generalization of theorem on Bohr for bases in spaces of holomorphic functions of several complex variables, \textit{J. Math. Anal.Appl.} \textbf{258}(2001), 429-447.

\bibitem{aleman-2012} {\sc A. Aleman} and {\sc A. Constantin}, Harmonic maps and ideal fluid flows, {\it Arch. Ration. Mech. Anal.} {\bf 204} (2012), 479--513.

\bibitem{Ali-Abhulhadi-Ng-CVEE-2016} {\sc R. M. Ali, Z. Abdulhadi} and {\sc C. Ng}, The Bohr radius for starlike loghramonic mappings, \textit{Complex Var. Elliptic Equ.} \textbf{61}(1) (2016), 1-14.

\bibitem{Ali-Ng-CVEE-2018} {\sc R. M. Ali} and {\sc Z. C. Ng}, The Bohr inequality in the hyperbolic plane, \textit{Complex Var. Elliptic Equ.} \textbf{63}(11)(2018), 1539-1557.

\bibitem{Ali-Abdulhadi-Ng-CVEE-2016} {\sc R. M. Ali, Z. Abdulhadi} and {\sc Z. C. Ng}, The Bohr radius for starlike logharmonic mappings, \textit{Complex Var. Elliptic Equ.} \textbf{61}(1)(2016), 1-14.

\bibitem{Alkhaleefah-Kayumov-Ponnusamy-PAMS-2019} {\sc S. A. Alkhaleefah, I. R. Kayumov} and {\sc S. Ponnusamy}, On the Bohr inequality with a fixed zero coefficients, \textit{Proc. Amer. Math. Soc.} \textbf{147}(12) (2019), 5263-5274.

\bibitem{Dahlnere-Khavinshon-CMFT-2004} {\sc C. B$ \acute{e}n\acute{e} $teau, A. Dahlner} and {\sc D. Khavinson,} Remarks on the Bohr phenomenon, \textit{Comput. Methods Funct. Theory}, \textbf{4}(2004), 1-19.

\bibitem{Aytuna-Djkov-BLMS-2013} {\sc A. Aytuna} and {\sc P. Djakov}, Bohr property of bases in the space of entire functions and its generalizations, \textit{Bull. London Math. Soc.} \textbf{45}(2)(2013), 411-420.

\bibitem{Bhowmik-Das-JMAA-2018} {\sc B. Bhowmik} and {\sc N. Das,} Bohr phenomenon for subordinating families of certain univalent functions, J. Math. Anal. Appl. 462 (2018), 1087-1098.

\bibitem{Bhowmik-Das-2019} {\sc B. Bhowmik} and {\sc N. Das,} Bohr phenomenon for operator valued functions with fixed initial coefficients, https://arxiv.org/pdf/2003.05810.pdf.

\bibitem{Bhowmik-Das-arxive-2020} {\sc B. Bhowmik} and {\sc N. Das,} On the Bohr phenomenon for complex valued and vector valued functions, https://arxiv.org/pdf/2011.12766.pdf.

\bibitem{Blasco-2009} {\sc O. Blasco}, The Bohr radius of a Banach space. In: Curbera, G.P., Mockenhaupt, G., Ricker,
W.J (eds.) Vector measures, Integration and Related Topics. Opearator theory: advances and applications, vol. 201, pp-59-64. Birkhauser Verlag, Basel (2009).

\bibitem{Boas-Khavinshon-PAMS-1997} {\sc H. P. Boas} and {\sc D. Khavinson}, Bohr's power series theorem in several variables, \textit{Proc. Amer. Math. Soc.} \textbf{125} (1997), 2975-2979.

\bibitem{Bohr-PLMS-1914} {\sc H. Bohr}, A theorem concerning power series, \textit{Proc. Lond. Math. Soc.} s2-13 (1914), 1-5.


\bibitem{constantin-2017} {\sc A. Constantin} and {\sc M. J. Martin}, A harmonic maps approach to fluid flows, {\it Math. Ann.} {\bf 369} (2017), 1--16.

\bibitem{Defant-2003} {\sc A. Defant, D. Garcia} and {\sc M. Maestre}, Bohr's power series theorem and local Banach space theory, \textit{J. Reine Angew. Math.} \textbf{557} (2003), 173-197.

\bibitem{Dixon-BLMS-1995} {\sc P. G. Dixon}, Banach algebras satisfying the non-unital von Neumann inequality, Bull. London Math. Soc. \textbf{27}(4)(1995), 359-362.

\bibitem{Duman-2011} {\sc E. Y. Duman}, Some distortion theorems for starlike logharmonic functions, RIMS Kokyuroku,
2011: 1772:1-7.

\bibitem{Evdoiridis-Ponnusamy-IM-2019} {\sc S. Evdoridis} and {\sc S. Ponnusamy}, Improved Bohr's inequality for locally univalent harmonic mappings, \textit{Indag. Math.} (N.S.) \textbf{30} (2019), 201-213.

\bibitem{Fournier-JMAA-2008} {\sc R. Fournier}, Asymptotics of the Bohr radius for polynomials of fixed degree, \textit{J. Math. Anal. Appl.} \textbf{338}(2008), 1100-1107.

\bibitem{Gong-1991} {\sc S. Gong}, Bieberbach Conjecture, Studies in Advanced Mathematics, \textit{Amer. Math. Soc.}, Providence, RI, 1991.

\bibitem{Jun-PAMS-1993}  {\sc S. H. Jun}, Univalent harmonic mappings on $ \Delta=\{z: |z|>1\} $, \textit{Proc. Amer. Math. Soc.} \textbf{119}(1)(1993), 109–114.

\bibitem{Kayumov-Ponnusamy-MN-2018} {\sc I. R. Kayumov, S. Ponnusamy} and {\sc N. Shakirov}, Bohr radius for localy univalent harmonic mappings, \textit{Math. Nachr.} \textbf{291} (2018), 1757-1768.

\bibitem{Kayumov-Ponnusamy-JMAA-2018} {\sc I. R. Kayumov} and {\sc S. Ponnusamy}, Bohr's inequalities for the analytic functions with lacunary series and harmonic functions, \textit{J. Math. Anal. Appl.} \textbf{465} (2018), 857-871.

\bibitem{Kayumov-Ponnusamy-CAMS-2020} {\sc I. R. Kayumov} and {\sc S. Ponnusamy}, Improved version of Bohr's inequalities, \textit{C. R. Math. Acad. Sci. Paris}, \textbf{358}(5) (2020), 615-620.

\bibitem{Kayumov-Ponnusamy-Shakirov-MN-2018} {\sc I. R. Kayumov, S. Ponnusamy} and {\sc N. Shakirov}, Bohr radius for locally univalent harmonic mappings, \textit{Math. Nachr}, \textbf{291}(2018), 1757-1768.

\bibitem{Li-Ponnusamy-Wang-BMMS-2013} {\sc P. Li, S.Ponnusamy} and {\sc X. Wang}, Some Properties of Planar p-Harmonic and log-p-Harmonic Mappings, \textit{Bull. Malyesian Math. Soc.} \textbf{36}(3)(2013), 595-609.

\bibitem{Liu-Ponnusamy-2020} {\sc Z. H. Liu} and {\sc S. Ponnusamy}, Some properties of univalent log-harmonic mappings, arXiv preprint arXiv:1808.07393.

\bibitem{Liu-Ponnuamy-arxive-2019} {\sc Z. H. Liu} and {\sc S. Ponnusamy}, On univalent log-harmonic mappings, arXiv preprint arXiv:1905.10551.

\bibitem{Liu-Ponnusamy-PAMs-2020} {\sc M. S. Liu} and {\sc S.  Ponnusamy}, Multidimensional analogues of refined Bohr's inequality, \textit{Proc. Amer. Math. Soc.} (2020) (to appear).

\bibitem{Mao-Ponnusamy-Wang-CVEE-2013} {\sc Z. Mao, S Ponnusamy} and {\sc X Wang}, Schwarzian derivative and Landau's theorem for logharmonic mappings, Complex Var. Elliptic Equ. \textbf{58}(8)(2013), 1093-1107.

\bibitem{Paulson-Popescu-Singh-2002} {\sc V. I.Paulsen, G. Popescu} and {\sc D. Singh}, On Bohr's inequality, \textit{Proc. Lond. Math. Soc.} s3-85 (2002), 493-512.

\end{thebibliography}
\end{document}